\documentclass[article, 12pt]{amsart}
\usepackage{latexsym}
\usepackage{amssymb, amsmath, amsthm}
\usepackage{graphicx}
\usepackage{subfigure}
\usepackage{tikz}

\newtheorem{theorem}{Theorem}[section]
\newtheorem{lemma}[theorem]{Lemma}
\theoremstyle{remark}
\newtheorem{remark}[theorem]{\bf Remark}

\newtheorem{example}[theorem]{\bf Example}
\newtheorem{proposition}[theorem]{\bf Proposition}

\newtheorem{definition}[theorem]{\bf Definition}

\newcommand{\f}{\mathcal{F}}
\newcommand{\g}{\mathcal{G}}
\newcommand{\p}{\mathbb{P}}

\newcommand{\aut}{{\rm Aut}{(\p)}}
\newcommand{\dd}{\mathcal{D}}
\newcommand{\ov}{\overrightarrow}
\def\acts{\curvearrowright}
\newcommand{\iLim}{\varprojlim}

\begin{document}

\title[Large orbits and the pseudo-arc]{Large orbits, projective Fra\"{i}ss\'{e} limits, and the pseudo-arc}

\author{Aleksandra Kwiatkowska}


\begin{abstract}
We show that the conjugacy action of the automorphism group, $\aut$, of the projective Fra\"{i}ss\'{e} limit $\p$, whose natural quotient is the pseudo-arc, on the set of involutions of $\p$, has a comeager orbit. 
\end{abstract}

\maketitle

\section{Introduction}

\subsection{The pseudo-arc}

The {\em pseudo-arc} $P$ is the unique hereditarily indecomposable chainable continuum.
Recall that a {\em continuum} is a compact and connected metric space; it is {\em indecomposable} if it is not a
union of two proper subcontinua, and it is {\em hereditarily indecomposable} if every subcontinuum
is indecomposable. We call a continuum  {\em chainable} if each open cover of it is refined by an open cover $U_1,U_2,\ldots, U_n$  such that for $i,j $,
$U_i\cap U_j\neq\emptyset$  if and only if $|j-i|\leq 1$.

The pseudo-arc has a remarkably rich structure, for example, it is (injectively) homogeneous; see \cite{B} and \cite{Leh}.
Irwin and Solecki \cite{IS} discovered that it is also projectively ultrahomogeneous.
Moreover, the collection of all subcontinua of $[0,1]^\mathbb{N}$  homeomorphic to the pseudo-arc is comeager in the space of all subcontinua of $[0,1]^\mathbb{N}$, 
equipped with the Hausdorff metric.
For more information on the pseudo-arc, see \cite{Lew}.

\subsection{Projective Fra\"{i}ss\'{e} theory}\label{sone}

We recall here basic notions and results on the projective Fra\"{i}ss\'{e} theory,   developed by Irwin and Solecki in \cite{IS}.

 Given a language $L$ that consists of relation symbols $\{r_i\}_{i\in I}$, and function symbols $\{f_j\}_{\in J}$,
a \emph{topological $L$-structure} is a compact zero-dimensional second-countable space~$A$ equipped with
closed relations $r_i^A$ and continuous functions $f_j^A$, $i\in I, j\in J$.
A~continuous surjection $\phi\colon B\to A$ is an
 \emph{epimorphism} if it preserves the structure, more precisely, for a function symbol $f$ of arity $n$ and $x_1,\ldots,x_n\in B$ we require:
\[
 f^A(\phi(x_1),\ldots,\phi(x_n))=\phi(f^B(x_1,\ldots,x_n));
\]
and for a relation symbol $r$ of arity $m$ and $x_1,\ldots,x_m\in B$  we require:
\begin{equation*}
\begin{split}
&  r^A(x_1,\ldots,x_m) \\ 
&\iff \exists y_1,\ldots,y_m\in B\left(\phi(y_1)=x_1,\ldots,\phi(y_m)=x_m, \mbox{ and } r^B(y_1,\ldots,y_m)\right).
\end{split}
\end{equation*}
By an \emph{isomorphism}  we mean a bijective epimorphism.

For the rest of this section fix a language $L$.
Let $\mathcal{F}$ be a family of finite topological\\ $L$-structures. We say that $\mathcal{F}$ is a \emph{projective Fra\"{i}ss\'{e} family}
if the following two conditions hold:

(F1) (the joint projection property: JPP) for any $A,B\in\f$ there are $C\in \f$ and epimorphisms from $C$ onto $A$ and from $C$ onto $B$;

(F2) (the amalgamation property: AP) for $A,B_1,B_2\in\f$ and any epimorphisms $\phi_1\colon B_1\to A$ and $\phi_2\colon B_2\to A$, there exist $C\in\f$,
 $\phi_3\colon C\to B_1$, and $\phi_4\colon C\to B_2$ such that $\phi_1\circ \phi_3=\phi_2\circ \phi_4$.

A topological $L$-structure $\mathbb{P}$ is a \emph{projective Fra\"{i}ss\'{e} limit} of $\mathcal{F}$ if the following three conditions hold:

(L1) (the projective universality) for any $A\in\f$ there is an epimorphism from $\mathbb{P}$ onto~$A$;

(L2) for any finite discrete topological space $X$ and any continuous function
 $f\colon \mathbb{P} \to X$ there are $A\in\f$, an epimorphism $\phi\colon \mathbb{P}\to A$, and a function
$f_0\colon A\to X$ such that $f = f_0\circ \phi$.

(L3) (the projective ultrahomogeneity) for any $A\in \f$ and any epimorphisms $\phi_1\colon \mathbb{P}\to A$ and $\phi_2\colon \mathbb{P}\to A$
there exists an isomorphism $\psi\colon \mathbb{P}\to \mathbb{P}$ such that $\phi_2=\phi_1\circ \psi$;

Here is the fundamental result in the projective Fra\"{i}ss\'{e} theory: 
\begin{theorem}[Irwin-Solecki, \cite{IS}]
 Let $\f$ be a countable projective Fra\"{i}ss\'{e} family of finite topological $L$-structures. Then:
\begin{enumerate}
 \item there exists a projective Fra\"{i}ss\'{e} limit of $\f$;\\
\item any two topological $L$-structures that are projective Fra\"{i}ss\'{e} limits are isomorphic.
\end{enumerate}
\end{theorem}

In the proposition below we state some properties of the projective Fra\"{i}ss\'{e} limit.

\begin{proposition}\label{fraisse}
\begin{enumerate}
 \item If $\mathbb{P}$ is the projective Fra\"{i}ss\'{e} limit, then the following condition (called \emph{the extension property}) holds:
Given $\phi_1\colon B\to A$, $A,B\in\f$, and $\phi_2\colon \mathbb{P}\to A$, then, there is $\psi\colon \mathbb{P}\to B$ such that $\phi_2=\phi_1\circ \psi$.\\
\item If $\mathbb{P}$ satisfies  the projective universality (L1), the extension property, and (L2), then it also satisfies  projective ultrahomogeneity, and therefore
is isomorphic to the projective Fra\"{i}ss\'{e} limit.
\end{enumerate}
\end{proposition}

For a sequence $\{A_n\}$ of finite structures and for epimorphisms $\alpha_n\colon A_{n+1} \to A_n$ we denote the inverse sequence by $\{A_n,\alpha_n\}$. For $m > n$ we let $\alpha_n^m= \alpha_n\circ \dots \circ \alpha_{m-1}$ and note that $\alpha_n^{n+1} = \alpha_n$. Given an inverse sequence $\{A_n,\alpha_n\}$, the associated inverse limit space is the subspace of the 
product space $\Pi A_i$ equal to the $\{(x_1,x_2,\ldots) : x_i \in A_i \text{ and } x_i = \alpha_i(x_{i+1})\}$, and is denoted as 
$\mathbb{P}=\iLim\{A_n,\alpha_n\}$. 
 We denote the canonical projection from the inverse limit space $\mathbb{P}$ to the $n$th factor space $A_n$ by $\alpha_n^{\infty}$.

We will use   \cite[Proposition~2.3]{B-C}, which can also be derived from the proof of \cite[Theorem 2.4]{IS}. 
\begin{proposition}\label{B-CThm}
Let $\f$ be a projective Fra\"{\i}ss\'e class. Let $\{A_n,\alpha_n\}$ be an inverse sequence in 
$\f$. Assume that for each $A \in \f$, $n \in \mathbb N$, and an epimorphism 
$f\colon A \to A_n$, there exist $m \geq n$ and an epimorphism $g\colon A_m \to A$ such that 
$f\circ g=\alpha^m_n$. Then $\iLim\{A_n,\alpha_n\}$ is the projective Fra\"{\i}ss\'e limit of $\f$.
\end{proposition}

\subsection{The pseudo-arc as a  projective Fra\"{i}ss\'{e} limit}

Let $L_0$ be the language that consists of one binary relation symbol $r$.
Let $\g$ denote the family of finite reflexive linear graphs, more precisely, we say that $A=([n],r^A)$, where $[n]=\{1,2,\ldots,n\}$ is a finite reflexive
 linear graph if  $r^A(x,y)$ holds if and only if
$x=y$, or $x=i,y=i+1$ for some $i=1,2,\ldots, n-1$, or  $x=i+1,y=i$ for some $i=1,2,\ldots, n-1$.

Recall the following results obtained by Irwin and Solecki \cite{IS}.
\begin{theorem}[Theorem 3.1 in \cite{IS}]
The family $\g$ is a projective Fra\"{i}ss\'{e} family.
\end{theorem}

\begin{lemma}[Lemma 4.1 in \cite{IS}] Let $\p$ be the projective Fra\"{i}ss\'{e} limit of $\g$. Then $r^{\mathbb{P}}$ is an equivalence
relation  whose each equivalence class has at most two elements.
\end{lemma}

\begin{theorem}[Theorem 4.2 in \cite{IS}]
Let $\p$ be the projective Fra\"{i}ss\'{e} limit of $\g$. Then $\p/r^{\p}$ is the pseudo-arc.
\end{theorem}

\subsection{Results}
A function $f\colon\p\to\p$ is an automorphism if and only if it is a homeomorphism and for every $x,y\in\p$, $r^\p(x,y) \iff r^\p(f(x),f(y))$. We denote the group of automorphisms of $\p$ by $\aut$. The topology on $\aut$ is inherited from the compact-open topology on the homeomorphism group $H(\p)$. 
Let 
\[ I(\p)=\{f\in \aut\colon f^2=\mbox{id}\}\]
be the set of involutions of $\p$. The group $\aut$ acts on $I(\p)$ by conjugation  $h\mapsto ghg^{-1}$, $g\in\aut$, and we will 
denote this action 
by $\aut\acts I(\p)$. Note that $I(\p)$ is a closed invariant subspace 
for the conjugacy action of $\aut$ on $\aut$. 
The main result of the present paper is the following theorem.
\begin{theorem}\label{abcd}
The action $\aut\acts I(\p)$ has a comeager orbit.
\end{theorem}
\noindent The proof of Theorem \ref{abcd} will be given in Section 3. 
In Section 2,  we show that $\aut\acts I(\p)$ has a dense orbit as well as we see how the existence  of a dense orbit of  $\aut\acts I(\p)$ implies the existence of  a dense orbit of the corresponding action of  the homeomorphism group of the pseudo-arc $P$ on the set of involutions of $P$.
In Appendix A, we present criteria for the automorphism group of a projective Fra\"{i}ss\'{e} limit and its conjugacy invariant closed subspace to have, respectively, a dense orbit and a comeager orbit 
in the conjugacy action.

\section{Dense orbits in the conjugacy actions of $\aut$ and $H(P)$}

Let $P$ denote the pseudo-arc and let $H(P)$ denote the homeomorphisms group of $P$, equipped with the compact-open topology. Let $I(P)$ be the set of involutions  of $P$.
The group $H(P)$ acts by conjugation on $I(P)$, and we denote this action by $H(P)\acts I(P)$.
In this section we prove the following theorems.
\begin{theorem}\label{pseu}
The action $H(P)\acts I(P)$ has a dense orbit.
\end{theorem}

\begin{theorem}\label{aut}
The action $\aut\acts I(\p)$ has a dense orbit.
\end{theorem}

\begin{remark}
We do not know if 
 $\aut$, has a dense conjugacy class or if $H(P)$ has a dense conjugacy class. All  existing "proofs" seem to have a flow. 
\end{remark}

Let $L=L_0\cup \{s\}=\{r,s\}$, where $s$ is a symbol for a binary relation.
 With some abuse of notation, we will be writing $(A,s^A)$, where $A=(A,r^A)$ is a topological $L_0$-structure, instead of $(A,r^A,s^A)$, whenever $(A,r^A,s^A)$ is a 
 topological $L$-structure.

Let $f\in\aut$. We view $(\p,f)$ as a topological $L$-structure $(\p,s^\p)$, where 
$s^\p(x,y) \iff (x,y)\in\mbox{graph}(f) \iff f(x)=y$.
Define \begin{equation*}
\begin{split}
\f= & \{ (A, s^A)\colon A\in\g \mbox{ and } \exists\phi\colon\p\to A \exists f\in I(\p) \mbox{ such that}\\
 &\phi\colon (\p,f)\to(A,s^A) \mbox{ is an epimorphism}\}.
\end{split}
\end{equation*}

\begin{remark}
Let $f\in \aut$. Let $A\in\g$. For a given epimorphism $\phi\colon\p\to A$ we can talk about a restriction of $f$ to $A$:
\[
 f\restriction A=\{(a,b)\in A^2\colon \phi^{-1}(a)\cap f(\phi^{-1}(b))\neq\emptyset\}.
\]
It is not difficult to see that
\[
 \f= \{ (A,s^A)\colon A\in\g \mbox{ and } \exists\phi\colon\p\to A \exists f\in I(\p)\ f\restriction A=s^A\}.
\]
\end{remark}

\begin{lemma}\label{replace}
Let $(A,s^A)\in\f$. Let $\phi\colon\p\to A$ be such that $\phi\colon (\p,f)\to(A,s^A) $ is an epimorphism. Let $\psi\colon \p\to A$ be an epimorphism.
Then there is $g\in\aut$ such that $\psi\colon (\p,gfg^{-1})\to (A,s^A)$ is an epimorphism.
\end{lemma} 
\begin{proof}
Using the projective universality, get $g\in\aut$  such that $\psi\circ g=\phi$.  This $g$ works.
\end{proof}
We will use several times Lemma \ref{replace} in  proofs of Propositions \ref{dense} and \ref{comeager} without mentioning it. We need some more results, which will be used later.

Let $A\in\g$ and let $s^A$ be a binary relation on $A$.
We say that $s^A$ is the {\em  antidiagonal} of $A=[n]=\{1,2,\ldots,n\}$ if $s^A=\{(k,n+1-k)\colon k=1,2,\ldots, n\}$. 
\begin{proposition}\label{lifttoinv}
Let $A\in\g$, let $s^A$ be the antidiagonal of $A$,  and let $\phi\colon \p\to A$ be an epimorphism. Then there is $h\in I(\p)$ such that $\phi\colon (\p,h)\to (A,s^A)$ is an epimorphism.
\end{proposition}
We first need a lemma.
\begin{lemma}\label{switch}
Let $A=[k],B=[m]\in\g$ and let $\phi\colon B\to A$ be an epimorphism. Then there is $D=[n]\in\g$ 
and  $\psi\colon D\to B$ such that for every $1\leq i\leq n$, we have
$k+1-\phi\circ\psi(i)=\phi\circ\psi(n+1-i)$.
\end{lemma}
\begin{proof} 
Let $\mbox{inv}_A\colon A\to A$ be given by $\mbox{inv}_A(i)=k+1-i$ for all $i=1,2,\ldots, k$.
Amalgamate $\phi$ and $\mbox{inv}_A\circ \phi$ and let $C\in\g$, and epimorphisms $\alpha, \beta\colon C\to  B$
 be such that $\phi\circ\alpha=\mbox{inv}_A\circ \phi\circ\beta$.

{\bf{ Claim.}} There is an epimorphism $\gamma\colon C'\to C$ such that, denoting $C'=[n']$,  
we have that $\alpha\circ\gamma(n')$ is adjacent to $\beta\circ\gamma(n')$.  

\begin{proof}[Proof of the Claim]
Take $n_1$ and an epimorphism $\gamma_1\colon C_1=[n_1]\to C=[n]$ such that $\gamma_1\restriction [n] = \mbox{id}$ and $\alpha\circ\gamma_1(n_1)=m$. Then take $n_2$ and an epimorphism $\gamma_2\colon C_2=[n_2]\to C_1=[n_1]$ such that $\gamma_2\restriction [n_1] = \mbox{id}$ and $\beta\circ\gamma_1\circ\gamma_2(n_2)=m$. We let $\gamma=\gamma_1\circ \gamma_2$.

Note that \[ \alpha\circ\gamma(n_1) - \beta\circ\gamma(n_1) = m -  \beta\circ\gamma(n_1) \geq 0\] and 
 \[ \alpha\circ\gamma(n_2) - \beta\circ\gamma(n_2) =  \alpha\circ\gamma(n_2) - m \leq 0.\]
Therefore there is $n_1\leq n'\leq n_2$ such that $\alpha\circ\gamma(n')$ is adjacent to $\beta\circ\gamma(n')$.
Take $C'=[n']$ and restrict $\gamma$ to $[n']$.
\end{proof}
Take $n=2n'$ and $D=[n]$. The $\psi\colon D\to B$ defined as follows:
\begin{equation*}
\psi(i) =
\begin{cases}
\alpha\circ\gamma (i) & \text{if } 1\leq i\leq n',\\
\beta\circ\gamma(n-i+1) & \text{if } n'+1\leq i\leq n,\\
  \end{cases}
\end{equation*}

is as required. Indeed, it is an epimorphism due to Claim. 
Furthermore, if $1\leq i\leq n'$, then
\begin{equation*}
\begin{split}
& k+1-\phi\circ\psi(i)=k+1-\phi\circ\alpha\circ\gamma(i)=k+1-\mbox{inv}_A\circ \phi\circ\beta\circ\gamma(i)\\
& =\phi\circ\beta\circ\gamma(i)
=\phi\circ\beta\circ\gamma(n - (n-i+1) +1) 
=\phi\circ\psi(n+1-i)
\end{split}
\end{equation*}
\end{proof}

\begin{proof}[Proof of Proposition \ref{lifttoinv}]
Thanks to lemma \ref{switch}, we can write $\p$ as the inverse limit of an inverse sequence $\{A_n,f_n\}$, where $f_n\colon A_{n+1}=[k_{n+1}]\to A_n=[k_n]$ is an epimorphism such that 
for every $1\leq i \leq k_{n+1}$, we have $k_n+1-f_n(i)=f_n(k_{n+1}+1-i)$, $A_1=A$, and $f^\infty_1=\phi$. Indeed, write first $\p$ as the inverse limit of an inverse sequence $\{B_n,g_n\}$ with a property that for any $m$ and an epimorphism $\alpha\colon C\to B_m$ there 
is $n>m$ and $\beta\colon B_n\to C$ such that $\alpha\circ\beta=g^n_m$, which we can do by Proposition \ref{B-CThm}. Then alternately apply this property of $\{B_n, g_n\}$ and Lemma \ref{switch} to construct the required $\{A_n,f_n\}$.

Once we have the inverse sequence $\{A_n,f_n\}$, take $h\colon\p\to \p$ given by
$h(a_1,a_2,a_3,\ldots )=(k_1+1-a_1, k_2+1-a_2, k_3+1-a_3,\ldots)$. This $h$ is well defined and is the required epimorphism.

\end{proof}

We split the proof of Theorem \ref{aut} into two propositions.

Recall from Introduction  that $\f$ has the JPP if and only if for every $(A,s^A),(B,s^B)\in \f$ there is $(C,s^C)\in\f$ and epimorphisms from  $(C,s^C)$ onto $(A,s^A)$
and from $(C,s^C)$ onto $(B,s^B)$.
\begin{proposition}\label{jpp}
 The family $\f$ has the JPP.
\end{proposition}
 
\begin{proposition}\label{dense}
 The property JPP for $\f$ implies that $\aut\acts I(\p)$ has a dense orbit.
\end{proposition}

 The proof of Proposition \ref{dense} will be an adaptation to our context of the proof of one of the directions of Theorem 2.1 in \cite{KR}.

 For $(A,s^A)\in\f$ and an epimorphism $\phi\colon\p\to A$ define
\[
 [\phi,s^A]=\{f\in I(\p)\colon \phi\colon (\p,f)\to(A,s^A) \mbox{ is an epimorphism}\}.
\]

\noindent Sets of the form $[\phi,s^A]$ are clopen in $I(\p)$.

\begin{lemma}\label{basis}
 The family of all sets $[\phi,s^A]$, where $(A,s^A)\in \f$, is a basis of the topology on $I(\p)$.
\end{lemma}

\begin{proof}
Take  $g\in I(\p)$,  $\epsilon>0$, and
 $U=\{f\in I(\p)\colon \forall x \ d(f(x),g(x))<\epsilon\}$
($d$ is  any metric on the underlying set of $\p$).
This is an open set. We want to find a clopen neighborhood of $g$ that is of the form $ [\phi,s^A]$ and is contained in $U$.
For this, take  an arbitrary partition $Q$ of $\p$ of mesh $<\epsilon$ and let
$R=\{ q_0\cap g^{-1}(q_1)\colon q_0,q_1\in Q\}$.
Let $A$ be a refinement of $R$ such that $A$ together with the relation $r^A$ inherited from $r^\p$ is in $\g$ (condition (L2) guarantees the existence of such $A$).
Let $\phi$ be the natural projection from $\p$ to $A$. By the choice of $A$, this is an epimorphism.
We let $s^A=\{(p,r)\in A^2\colon \exists_{x\in p}\exists_{y\in r} g(x)=y\}$.
Clearly $g\in [\phi,s^A]$ and $(A,s^A)\in \f$.
Take any $p\in A$, say $p\subseteq q_0\cap g^{-1}(q_1)$, $q_0,q_1\in R$. Then $g(p)\subseteq q_1$.
Now take any $f\in [\phi,s^A]$ and notice that $f(p)\subseteq q_1$. Since $\mbox{diam}(q_1)<\epsilon$
  and  $p\in A$ was arbitrary, we get $f\in U$.
\end{proof}

For $(A,s^A)\in\f$ and an epimorphism $\phi\colon\p\to A$ define
\[
 D(\phi,s^A)=\{f\in I(\p)\colon \exists g\in\aut\ gfg^{-1}\in [\phi,s^A]\}.
\]
This set is open in $I(\p)$.
\begin{lemma}
The set $D(\phi,s^A)$, where $(A,s^A)\in\f$, is  dense in $I(\p)$.
\end{lemma}

\begin{proof}
Fix $D(\phi,s^A)$ and take $[\psi,s^B]$. We show that $D(\phi,s^A)\cap [\psi,s^B]\neq\emptyset$. Since sets $[\psi,s^B]$ form a basis, this will finish the proof.
Using the JPP, take $(C,s^C)\in\f$ and epimorphisms $\alpha\colon (C,s^C)\to (A,s^A)$ and $\beta\colon (C,s^C)\to (B,s^B)$.
Using the extension property, find $\gamma\colon \p\to C$ and $\delta\colon\p\to C$ such that  $\phi= \alpha\circ\gamma$ and
$\psi=\beta\circ \delta$. 
Note that $[\delta,s^C]\neq\emptyset$. Indeed, since $(C,s^C)\in\f$, there is $f_0\in I(\p)$ and an epimorphism $\phi_0\colon (\p,f_0)\to  (C,s^C)$. Hence, by projective ultrahomogeneity, for some $g_0\in\aut$, we have $g_0f_0g_0^{-1}\in [\delta,s^C]$.
Take any $f\in[\delta,s^C]\subseteq [\psi,s^B]$ and take $g\in\aut$ such that $gfg^{-1}\in [\gamma,s^C]$.
Then   $gfg^{-1}\in [\phi,s^A]$, and therefore $D(\phi,s^A)\cap [\psi,s^B]\neq\emptyset$. 
\end{proof}

\begin{proof}[Proof of Proposition \ref{dense}]
The intersection of all $D(\phi,s^A)$ is open and dense, in particular it is nonempty. From the definition of $D(\phi,s^A)$,
every automorphism in this intersection has a dense orbit.
\end{proof}

\begin{definition}
Let $A\in\g$ or $A=\p$. Let $s^A$ be a binary relation on $A$.

(1) We say that $s^A$ is {\em surjective} if for every $a\in A$ there are $b,c\in A$ such that $s^A(a,b)$ and $s^A(c,a)$.

(2) We say that $s^A$ is {\em connected} if the graph $G_{(A,s^A)}$ that has $s^A=\{(a,b)\in A^2\colon s^A(a,b)\}$ as the set of vertices and $\{((a,b),(c,d))\colon r^A(a,c), r^A(b,d)\}$ as the
 set of edges, satisfies the following: for any clopen $X\subseteq s^A$ 
such that $X$ and $s^A\setminus X$ are nonempty there are  $x\in X$ and $y\in s^A\setminus X$
such that $x$ and $y$ are joined by an edge.
(The topology on $s^A$ is inherited from the product topology on $A\times A$. If $A\in\g$, then $A\times A$ is discrete.) 

(3) We say that $s^A$ is {\em symmetric} if for every $x,y\in A$, it holds $(x,y)\in s^A$ iff $(y,x)\in s^A$.
\end{definition}

\begin{remark}
Note that for  $A\in \g$, $s^A$ is connected if and only if $G_{(A,s^A)}$ is connected as a graph (that is, every two vertices are connected by a path).
\end{remark}

\begin{example}
(1) Take $A=(\{1,2,3,4\},s^A)$, where $s^A=\{(1,3),(2,3),(3,1),(3,2),(3,4),\\ (4,1)\}$. Then $s^A$ is connected.

(2) Take $A=(\{1,2,3,4\},s^A)$, where $s^A=\{(1,2),(2,1),(2,4),(3,3),(3,4),(4,2)\}$. Then $s^A$ is not connected.
\end{example}

To prove Proposition \ref{jpp} we need Lemmas \ref{slawek} and \ref{pse}.
The next lemma is due to  Solecki.
\begin{lemma}\label{slawek}
 For any  $(A,s^A)$ with $A\in\g$ and $s^A$ surjective, symmetric, and connected there is $(B,s^B)$ such that: $s^B$ is the antidiagonal of $B$ and there exists an epimorphism from $(B,s^B)$ onto $(A,s^A)$ .
\end{lemma}

\begin{proof}[Proof of Lemma \ref{slawek}]
Take  $(A,s^A)$ with $A\in\g$ and $s^A$ surjective, symmetric, and connected. Write $A=[k]$.
 We let $[k]\times [k]$ to be the product graph.
 Since $s^A$ surjective and connected, there is $ i_0$ such that $ (i_0, i_0)\in s^A$  or  $(i_0, i_0+1)\in s^A$
 or  $(i_0+1, i_0)\in s^A$. 
Let $ h\colon [m]\to s^A$, for some m, be a surjective graph homomorphism with
 $h(1) = (i_0, i_0)$ or  $h(1) = (i_0, i_0+1)$
 or $h(1) = (i_0+1, i_0)$, according to the cases above.
Define $\phi\colon [4m]\to [k]$ as follows: 

\begin{equation*}
\phi(i) =
\begin{cases}
\pi_1(h(i)) & \text{if } 1\leq i\leq m,\\
\pi_1(h(2m-i+1)) & \text{if } m+1\leq i\leq 2m,\\
  \pi_2(h(i-2m)) & \text{if } 2m+1\leq i\leq 3m,\\
\pi_2(h(4m-i+1)) & \text{if } 3m+1\leq i\leq 4m.
\end{cases}
\end{equation*}

 Let $ s^B\subseteq [4m]\times [4m]$ be the antidiagonal of $B=[4m]$.
Then, just by applying the above formulas, we see that $\phi\colon  (B,s^B)\to (A,s^A)$ is an epimorphism.

\end{proof}

\begin{lemma}\label{pse}
 Let $A\in\g$. Then $(A,s^A)\in\f $ if and only if
$s^A$ is surjective, symmetric, and connected.
\end{lemma}

\begin{proof}
Suppose that $(A,s^A)\in\f $.
Let $\phi\colon (\p,f)\to(A,s^A) $ be an epimorphism, where $f\in I(\p)$.
Since $f$ is a bijection, ${\rm graph}(f)=\{(x,f(x))\colon x\in\p\}$ is a surjective relation on $\p$. As $f\in I(\p)$, ${\rm graph}(f)$ is symmetric.
We show that $G_{(\p,{\rm graph}(f))}$ is connected. Suppose towards a contradiction that there is a clopen set $X\subseteq {\rm{ graph}(f))}$ such that $X$
and  ${\rm graph}(f)\setminus X$ are nonempty, and there are no $(x,f(x))\in X$, $(y,f(y))\in {\rm graph }(f)\setminus X$ such that 
 $r^{\p}(x,y)$ (and $r^{\p} (f(x),f(y))$). Let $Y$ be the projection of $X$ into the first coordinate. Then $ Y$ and $\p\setminus Y$ are nonempty clopen and for no 
 $x\in Y$ and $y\in \p\setminus Y$, $r^\p(x,y)$.
  However, this is impossible (apply (L2) to $A=\{Y, \p\setminus Y\}$ and the natural projection from $\p$ to $A$).
Finally, observe that surjective, symmetric, connected relations  are preserved by epimorphisms.

For the other direction, take  $(A,s^A)$ such that $A\in\g$ and $s^A$ is surjective, symmetric, and connected. Apply  Lemma \ref{slawek} and get $B$ and   an epimorphism $\psi\colon   (B,s^B)\to(A,s^A)$, where  $s^B$ is the antidiagonal of $B$.
Let $\phi\colon \p\to B$ be any epimorphism and 
apply Proposition \ref{lifttoinv} to $\phi$. Then $\psi\circ \phi$ witnesses that
$(A,s^A)\in\f $.

\end{proof}

\begin{proof}[Proof of Proposition \ref{jpp}]
 Take $(A,s^A),(B,s^B)\in\f$. Then $s^A$ and $s^B$ are surjective, symmetric, and connected. Without loss of generality, by Lemma \ref{slawek}, $s^A$ and $s^B$ are antidiagonals of $A$ and $B$, respectively. 
Write $A=[k]$ and $B=[n]$. Take $C=[kn]$ and let $s^C$ be the antidiagonal of $C$.
By Proposition \ref{lifttoinv}, $(C,s^C)\in\f$.
We show that this $(C,s^C)$ works.
For this, take $\phi_1\colon (C,s^C)\to (A,s^A)$ given by
$\phi_1((i-1)n+j)=i$, $i=1,2,\ldots, k$, $j=1,2,\ldots,n$, and take
$\phi_2\colon (C,s^C)\to (B,s^B)$ given by
$\phi_2((i-1)k+j)=i$, $i=1,2,\ldots, n$, $j=1,2,\ldots,k$.

\end{proof}

\begin{proof}[Proof of Theorem \ref{aut}]
Apply Propositions \ref{jpp} and \ref{dense}.
\end{proof}

Let $\pi\colon \p\to P$ be the canonical projection and let   $f\mapsto f^*$
 be the canonical embedding of $\aut$ to $H(P)$, see  page 3090 of \cite{IS}.
As proved in \cite[Lemma 4.8]{IS}, the image of this embedding is dense in $H(P)$.

We need one more lemma to prove Theorem \ref{pseu}.
 \begin{lemma}\label{denseI(p)}
The  $f\mapsto f^*$ is an embedding of $I(\p)$  into $I(P)$, which has a dense image.
 \end{lemma}
\begin{proof}
Since $f\mapsto f^*$ is an embedding of $\aut$  into $H(P)$, 
its restriction to $I(\p)$ is still an embedding. We have to show that its image is dense.

Pick $h\in I(P)$, let $\epsilon>0$, and let $d$ be a metric on $P$. The sup metric corresponding to $d$ we denote by $d_{sup}$. We will find $f\in I(\p)$ such that $d_{sup}(f^*,h)<\epsilon$. Let $0<\delta<\epsilon$ be such that for any $x,y\in P$ if $d(x,y)<\delta/2$ then $d(h(x), h(y))<\epsilon/2$. Pick $E\in\g$ and an epimorphism $\phi\colon \p\to E$  such that each set in the closed cover $\alpha=\{\pi(\phi^{-1}(e))\colon e\in E\}$ of $P$ has diameter $<\delta/2$. Take $s^{\alpha}=\{(A,B)\in\alpha^2\colon h(A)\cap B\neq\emptyset\}$
and let $s^E(e_1,e_2)$ iff  $s^{\alpha}(\pi(\phi^{-1}(e_1), \pi(\phi^{-1}(e_2)))$. Clearly the binary relation $s^E$ is surjective and symmetric. We now show that it is also connected. Since $P$ is connected, so is the graph of $h$, i.e. $G=\{(x,y)\colon h(x)=y\}$. Suppose that $\mathcal{C},\mathcal{D}$
 partition $s^\alpha$ and witness that  $s^E$ is not connected. Therefore  $\mathcal{C},\mathcal{D}\neq\emptyset$ and for any $(A_1,B_1)\in\mathcal{C}$ and  $(A_2,B_2)\in\mathcal{D}$ we have 
 $(A_1\times B_1)\cap (A_2\times B_2)=\emptyset $.
Then
\[C=\{(x,y)\in G\colon \mbox{for some } (A,B)\in  \mathcal{C} \mbox{ we have } x\in A, y\in B\}\] and  
\[D=\{(x,y)\in G\colon \mbox{for some } (A,B)\in  \mathcal{D} \mbox{ we have } x\in A, y\in B\}\] 
are closed nonempty sets that partition $G$.
This gives a contradiction and  completes the proof that $s^E$ is connected.

Let $\psi\colon F\to E$ be an epimorphism such that $\psi\colon (F, s^F)\to (E, s^E)$ is an epimorphism, where $s^F$ is the antidiagonal of $F$.  Let $p\colon\p\to F$ be such that 
$\psi\circ p=\phi$. 
Using Proposition \ref{lifttoinv}, take $f\in I(\p)$ such that $p\colon (\p,f)\to (F,s^F)$ is an epimorphism.
Then $f$ is as required.
Indeed, take $x\in P$ and $A\in\alpha$ such that $x\in A$ (there is either one or two such $A$'s). Let $B\in\alpha$ be such that $h(x)\in B$, hence 
$(A,B)\in s^\alpha$. Let $C\in\alpha$ be such that $f^*(x)\in C$. Since 
$\phi\colon (\p,f)\to (E,s^E)$ is an epimorphism, we get that $(A,C)\in s^{\alpha}$.
Take $x'\in A$ such that $h(x')\in C$. As $d(x,x')<\delta/2$, we get $d(h(x),h(x'))<\epsilon/2$. Therefore $d(h(x),f^*(x))\leq d(h(x),h(x'))+d(h(x'),f^*(x))\leq \epsilon/2 +\delta/2\leq \epsilon$.

\end{proof}

\begin{proof}[Proof of Theorem \ref{pseu}]
Let $f\in I(\p)$  have a dense orbit in the action $\aut\acts I(\p)$. Then $f^*$ has 
a dense orbit in the conjugacy action of $H(P)$ on $\overline{\pi(I(\p))}$.
By Lemma \ref{denseI(p)}, in fact, $f^*$ has a dense orbit in $I(P)$. 
\end{proof}

\section{Comeager orbit in $\aut\acts I(\p)$}

In  this section we show our main theorem.
\begin{theorem}\label{maint}
The conjugacy action $\aut\acts I(\p)$ has a comeager orbit.
\end{theorem}

We say   that $\f$ has the {\em coinitial  amalgamation property} (the CAP) if and only if for every $(A_0,s^{A_0})\in\f$ there is $(A,s^A)\in \f$ and an epimorphisms
$\psi\colon (A,s^A)\to (A_0,s^{A_0})$ such that for every $(B,s^B),(C,s^C)\in\f$ and epimorphism $\phi_1\colon (B,s^B)\to (A,s^A)$ and $\phi_2\colon (C,s^C)\to (A,s^A)$
there is $(D,s^D)\in\f$ and epimorphisms $\phi_3\colon (D,s^D)\to (B,s^B)$ and $\phi_4\colon (D,s^D)\to (C,s^C)$ such that $\phi_1\circ\phi_3=\phi_2\circ\phi_4$. 

We split the proof of Theorem \ref{maint} into Propositions \ref{ap} and \ref{comeager}.
\begin{proposition}\label{ap}
The family $\f$ has the CAP. More precisely, the coinitial family $\dd$, defined in Lemma \ref{ddd} below, has the AP.
\end{proposition}

Since we already know that $\f$ has the JPP, Proposition \ref{ap} together with Proposition~\ref{comeager} will finish the proof of Theorem \ref{maint}.
\begin{proposition}\label{comeager}
Properties  CAP and JPP for $\f$ imply that $\aut\acts I(\p)$ has a comeager orbit.
\end{proposition}

It will be convenient for us to work only with those structures   that have an even number of elements.
\begin{lemma}\label{ddd}
The family 
\[
\dd=\{(A,s^A)\in\f\colon s^A \mbox{ is the antidiagonal of $A$ and } |A| \mbox{ is an even number} \}
\]
 is coinitial in $\f$.
\end{lemma}

\begin{proof}
We know already that $\dd_0=\{(A,s^A)\in\f\colon s^A \mbox{ is the antidiagonal  of } A\}$ is coinitial in $\f$ (Lemma \ref{slawek}).
Take $(A,s^A)\in\dd_0$. Write $A=[k]$. Let  $B=[2k]$ and let $s^B$ be the antidiagonal of $B$. Take
$\phi\colon (B,s^B)\to (A,s^A)$ given by $\phi((i-1)n+j)=i$, $i=1,2,\ldots, k$, $j=1,2$. This is an epimorphism.
\end{proof}


 The proof of Proposition \ref{comeager} will be an adaptation to our context of the proof of one of the directions of Theorem 3.4 in \cite{KR}. In the proof we use the following proposition.

\begin{proposition}[Proposition 3.2 in \cite{KR}]\label{KR}
Let $G$ be a non-archimedean group acting continuously on a Polish space $X$ and let $x\in X$. Then the following are equivalent:
\begin{enumerate}
\item the  orbit $G\cdot x$ is non-meager;
\item for each open subgroup $V<G$, $V\cdot x$ is somewhere dense;
\item for each open subgroup $V<G$, $x\in \mbox{Int}(\overline{ V\cdot x})$.
\end{enumerate}
\end{proposition}

For $(A,s^A)\in\f$ and an epimorphism $\phi\colon\p\to A$ we say that $((B,s^B),\psi,\bar{\psi})$ is an {\em extension} of $((A,s^A),\phi)$ if 
$(B,s^B)\in\f$, $\psi\colon\p\to B$ is an epimorphism, $\bar{\psi}\colon (B,s^B)\to (A,s^A)$ is an epimorphism, and $\phi=\bar{\psi}\circ\psi$.

\begin{proof}[Proof of Proposition \ref{comeager}]
We show that there is $f\in I(\p)$ with a dense and non-meager orbit. Clearly such $f$ has a comeager orbit.
For $(A,s^A)\in\f$ and an epimorphism $\phi\colon\p\to A$ we defined in Section 2
\[
 [\phi,s^A]=\{f\in I(\p)\colon \phi\colon (\p,f)\to(A,s^A) \mbox{ is an epimorphism}\},
\]
\[
 D(\phi,s^A)=\{f\in I(\p)\colon \exists g\in\aut\ gfg^{-1}\in [\phi,s^A]\},
\]
and we showed that every $D(\phi,s^A)$ is open and dense in $I(\p)$.

We need some more definitions.
Let $(A,s^A)\in\f$ and  $\phi\colon\p\to A$ be an epimorphism. 
Let $\mbox{id}_A$ be the surjective relation on $A$ satisfying $\mbox{ id}_A(x,y) \iff x=y$. 
Let
\[
[\phi, \mbox{id}_A]_{\aut}=\{h\in\aut\colon \phi\colon (\p,h)\to(A,\mbox{id}_A) \mbox{ is an epimorphism}\}  
\]  and let for $f\in I(\p)$
\[
c(\phi,f)=\{gfg^{-1}\colon    
g\in[\phi, \mbox{id}_A]_{\aut}\}. 
\]

\noindent Let $((A_m,s^{A_m}),\phi_m, \bar{\phi}_m)$  list  all extensions of $((A,s^A),\phi)$ such that additionally \\ $(A_m,s^{A_m})\in\dd$. 
 Further, for a given $m$, let $((A_m^n,s^{A_m^n}),\phi_m^n, \bar{\phi}_m^n)$  list  all extensions of $((A_m,s^{A_m}),\phi_m)$.
Define

%
\[
 E(\phi,s^A)=  \{f\in I(\p)\colon \mbox{ if }  f\in [\phi,s^A]  
  \mbox{ then for some }  m,\  f\in [\phi_m,s^{A_m}]\}
\]

\noindent and let
\[
F_{m,n}(\phi,s^A)=  \{f\in I(\p)\colon \mbox{ if } f\in [\phi_m,s^{A_m}] 
  \mbox{ then } c(\phi,f)\cap [\phi_m^n,s^{A_m^n}]\neq\emptyset\}.
\]

 Observe that $E(\phi,s^A)$ and $F_{m,n}(\phi,s^A)$ are open, and by (L2), $E(\phi,s^A)$ is dense.
 Claims 1, 2, and 3 will finish the proof of the proposition.

\noindent {\bf{Claim 1.}}    
 The collection of all $[\phi_m,s^{A_m}]$ coming from extensions of  $((A,s^{A}),\phi)$ forms a basis of $[\phi,s^A]$.
\begin{proof}
 The proof of this claim goes along the lines of the proof of Lemma \ref{basis}.
\end{proof}
    
 \noindent   {\bf{Claim 2.}}
   The set $F_{m,n}(\phi,s^A)$ is dense.
   \begin{proof}
   Take $(B,s^B)\in\f$ and an epimorphism $\psi\colon\p\to B$. We show that $F_{m,n}(\phi,s^A)\cap[\psi,s^B]\neq\emptyset$.
    We can assume that for some $\psi$ and $\bar{\psi}$, $((B,s^{B}),\psi, \bar{\psi})$ is an extension of $((A_m,s^{A_m}),\phi_m)$.
    Since $(A_m,s^{A_m})\in\dd$, there is $(C,s^C)\in\f$ and there are $\alpha\colon (C,s^C)\to (A_m^n,s^{A_m^n})$ and
   $\beta\colon (C,s^C)\to (B,s^B)$ such that $\bar{\phi_m^n}\circ\alpha=\bar{\psi}\circ \beta$. Using the extension property for $\g$, take 
   $\gamma\colon\p\to C$ and $\delta\colon\p\to C$ such that $\alpha\circ\gamma=\phi^n_m$ and $\beta\circ\delta=\psi$.
   Take any $f\in [\delta, s^C]\subseteq [\psi,s^B]\subseteq [\phi_m,s^{A_m}]$.
   We want to show $f\in F_{m,n}(\phi,s^A)$.
   Take $g\in\aut$ such that $\gamma\circ g=\delta$. Then $gfg^{-1}\in [\gamma,s^C]\subseteq[\phi^n_m, s^{A_m^n}]$.
 Since $\bar{\phi}_m^n\circ\alpha\circ\gamma=\bar{\psi}\circ\beta\circ\delta=\phi_m$, we have
$g\in [\phi_m, \mbox{id}_{A_m}]_{\aut}\subseteq [\phi, \mbox{id}_{A}]_{\aut}$, so $gfg^{-1}\in c(\phi,f)  $, and we are done.
     \end{proof}
   
  \noindent  {\bf{Claim 3.}}
Whenever $f$ is in the intersection of all $D(\phi,s^A)$, $E(\phi,s^A)$, and $F_{m,n}(\phi,s^A)$, where     $(A,s^A)\in\f$ and  $\phi\colon\p\to A$ is an epimorphism,
then it has a comeager orbit.
\begin{proof}
We already know that such $f$ has a dense orbit. We show that the orbit of $f$ is also non-meager.
Since
$\{[\phi,\mbox{id}_{A}]_{\aut}\colon A\in\g, \ \phi\colon \p\to A \mbox{ is an epimorphism}\}$ forms a basis of the identity of $\aut$ that consists of open subgroups,
via Proposition \ref{KR},
it suffices to show that for a given $A\in\g$ and an epimorphism $\phi\colon \p\to A$, $c(\phi,f)$ is somewhere dense. 

Take $s^A$ satisfying $s^A(a,b)$ if and only if there are $x,y\in\p$ such that $\phi(x)=a$, $\phi(y)=b$, and $f(x)=y$. Then $f\in[\phi,s^A]$.
 Since $f\in E(\phi,s^A)$, for some $m$, $f\in[\phi_m,s^{A_m}]$.
 We claim that $c(\phi,f)$ is dense in $[\phi_m,s^{A_m}]$. This is because  $f\in F_{m,n}(\phi,s^A)$ implies $c(\phi,f)\cap [A_m^n,s^{A_m^n}]\neq\emptyset$, and  because sets
  $[\phi_m^n,s^{A_m^n}]$ form a basis  of $[\phi_m,s^{A_m}]$.
\end{proof}

\end{proof}


We will frequently denote structures in $\g$ by  $[m]=\{1,2,\ldots, m\}$, $[-k,-1]\cup[1,k]=\{-k,\ldots,-1,1,\ldots,k\}$, etc. 
From now on, whenever we write $s^A,s^B, s^C, s^D$, we  always mean the antidiagonal of $A,B,C,D$, respectively.

In the rest of the section we prove Proposition \ref{ap}. We illustrate our proof in Example \ref{exa}.
We start with a simple lemma. The proof is straightforward.
\begin{lemma}\label{drops}
Let $A=[-k,-1]\cup[1,k]$ and $B=[-l,-1]\cup[1,l]$. Let $\phi\colon B\to A$ be an epimorphism. Then $\phi\colon (B,s^B)\to (A,s^A)$ is an epimorphism
if and only if for every $i$, $\phi(i)=-\phi(-i)$. In particular, if $\phi\colon (B,s^B)\to (A,s^A)$  is an epimorphism, then $\phi(1)=1$ and $\phi(-1)=-1$, or
 $\phi(1)=-1$ and $\phi(-1)=1$. 
 \end{lemma}

Therefore, it is enough to show the following proposition.
\begin{proposition}\label{mainp}
Let $A=[-k,-1]\cup[1,k]$, $B=[-l,-1]\cup[1,l]$, and $C=[-m,-1]\cup[1,m]$. Let $\phi_1\colon B\to A$ and  $\phi_2\colon C\to A$ be epimorphisms
such that for every $i$, $\phi_1(i)=-\phi_1(-i)$ and $\phi_2(i)=-\phi_2(-i)$.
Then there  are $D=[1,n]$ and epimorphisms $\psi_1\colon D\to B$ 
and $\psi_2\colon D\to C$ such that $\phi_1\circ \psi_1=\phi_2\circ \psi_2$, and moreover $\psi_1(1)\in \{-1,1\}$ and $\psi_2(1)\in \{-1,1\}$.
\end{proposition}

\begin{proof}[Proof that Proposition  \ref{mainp} implies Proposition \ref{ap}]
Take $A=[-k,-1]\cup[1,k]$, $B=[-l,\\ -1]\cup[1,l]$, and $C=[-m,-1]\cup[1,m]$.  Take epimorphisms $\phi_1\colon (B,s^B)\to (A,s^A)$ and  $\phi_2\colon (C,s^C)\to (A,s^A)$.
From  Lemma \ref{drops}, for every $i$, $\phi_1(i)=-\phi_1(-i)$ and $\phi_2(i)=-\phi_2(-i)$. From the conclusion of Proposition \ref{mainp}  get $D=[1,n],\psi_1,\psi_2$.
Take $D'=[-n,-1]\cup[1,n]$ and let $s^{D'}$ be the antidiagonal of $D'$. Extend $\psi_1$ and $\psi_2$ to $D'$ so that 
for every $i$, $\psi_1(i)=-\psi_1(-i)$ and $\psi_2(i)=-\psi_2(-i)$. Then $\psi_1\colon (D', s^{D'})\to (B,s^B)$ and  $\psi_2\colon (D', s^{D'})\to (C,s^C)$
are epimorphisms and $\phi_1\circ \psi_1=\phi_2\circ \psi_2$.
\end{proof}

To show Proposition \ref{mainp}, we need  the Steinhaus' chessboard theorem. The Steinhaus' chessboard theorem was first used by Solecki to show the amalgamation
property of the family of finite reflexive linear graphs (see Remark \ref{remsla} for the sketch of his proof). We use the Steinhaus' chessboard theorem
as one of the ingredients of the proof of Proposition \ref{ap}.

For $m,n$ positive define a {\em chessboard} to be $\mathbb{C}=[m]\times [n]$. The 
{\em boundary} of the chessboard $\mathbb{C}$, denoted by
$\mbox{Bd}(\mathbb{C})$, is defined to be the set $\left(\{1,m\}\times[n]\right)\cup \left( [m]\times \{1,n\}\right)$. For $(a_1,b_1),(a_2,b_2)\in \mathbb{C}$, 
we say that they are {\em 8-adjacent} if they are different and $|a_1-a_2|\leq 1$, and $|b_1-b_2|\leq 1$; they are  {\em 4-adjacent} if they are different and either
 $|a_1-a_2|\leq 1$ and $b_1=b_2$, or $a_1=a_2$ and $|b_1-b_2|\leq 1$. 
A sequence $x_1,x_2,\ldots, x_l$ is a {\em 4-path} (an {\em 8-path})  
from $A\subseteq\mathbb{C}$ to $B\subseteq \mathbb{C}$ if $x_1\in A$, $x_l\in B$,
and for every $i$, $x_i$ and $x_{i+1}$ are   4-adjacent (8-adjacent).
For $x,y\in \mbox{Bd}(\mathbb{C})$, $x\neq y$, there are exactly two 4-paths from $x$ to $y$ such that every element in the path is in the boundary: clockwise and 
counter-clockwise. If $x=x_1,x_2,\ldots, x_l=y$ is the clockwise path from $x$ to $y$, we let $\ov{xy}=\{x_1,x_2,\ldots, x_l\}$.
For  $x\in \mbox{Bd}(\mathbb{C})$, we let $\ov{xx}=\{x\}$. We say that $w,x,y,z\in \mbox{Bd}(\mathbb{C})$ is an {\em oriented quadruple} if 
$y,z\notin \ov{wx}$ and $z\notin\ov{xy}$. A {\em coloring} is any function $f\colon \mathbb{C}\to \{\mbox{black,white}\}$.

The theorem below is due to Hugo Steinhaus, for the proof we refer the reader to \cite{S}. We use the  chessboard theorem to obtain various
amalgamation results (Lemma \ref{amalg}). 
\begin{theorem}[Steinhaus' chessboard theorem]
Let $\mathbb{C}$ be a chessboard. Let $w,x,y,z\in \mbox{Bd}(\mathbb{C})$  be an oriented quadruple. 
Then for every coloring $f\colon \mathbb{C}\to \{black,white\}$ the existence of an 8-path from $\ov{wx}$ to $\ov{yz}$ is equivalent to the non-existence of a
4-path from $\ov{xy}$ to $\ov{zw}$.
\end{theorem}

\begin{proposition}\label{path}
Let $A=[-k,-1]\cup[1,k]$, $B=[-l,-1]\cup[1,l]$, and $C=[-m,-1]\cup[1,m]$. Let $\phi_1\colon B\to A$ and  $\phi_2\colon C\to A$ be epimorphisms
such that for every $i$, $\phi_1(i)=-\phi_1(-i)$ and $\phi_2(i)=-\phi_2(-i)$.
We let $(i,j)\in B\times C$ to be black if and only if $\phi_1(i)=\phi_2(j)$.
Then there are: a black 8-path from $\{-1,1\}\times \{-1,1\}$ to $\{-l\}\times [-m,m]$,  a black 8-path from  $\{-1,1\}\times \{-1,1\}$ to $\{l\}\times [-m,m]$,
 a black 8-path from  $\{-1,1\}\times \{-1,1\}$ to $[-l,l]\times \{-m\}$, and a black 8-path from  $\{-1,1\}\times \{-1,1\}$ to  $[-l,l]\times \{m\}$. 
\end{proposition}
\begin{proof}[Proof that  Proposition \ref{path} implies Proposition \ref{mainp}]
Let $A=[-k,-1]\cup[1,k]$, $B=[-l,-1]\\ \cup[1,l]$, and $C=[-m,-1]\cup[1,m]$. Let $\phi_1\colon B\to A$ and  $\phi_2\colon C\to A$ be epimorphisms
such that for every $i$, $\phi_1(i)=-\phi_1(-i)$ and $\phi_2(i)=-\phi_2(-i)$.
Observe  that either $\{-1\}\times \{-1\}$ and $\{1\}\times\{1\}$ are black and $\{-1\}\times \{1\}$ and $\{1\}\times\{-1\}$ are white, or
$\{-1\}\times \{-1\}$ and $\{1\}\times\{1\}$ are white and $\{-1\}\times \{1\}$ and $\{1\}\times\{-1\}$ are black. Without loss of generality, the former holds.
Let $w_1,w_2,\ldots,w_p$ be a black 8-path from $\{1\}\times\{1\}$ to $\{-l\}\times [-m,m]$ (clearly, there is a black 8-path from $\{1\}\times\{1\}$ to $\{-l\}\times [-m,m]$
if and only if there is a black 8-path from $\{-1\}\times\{-1\}$ to $\{-l\}\times [-m,m]$), Let $x_1,x_2,\ldots,x_q$ be a black 8-path from $\{1\}\times\{1\}$ to 
$\{l\}\times [-m,m]$, let $y_1,y_2,\ldots,y_r$ be a black 8-path from  $\{1\}\times \{1\}$ to $[-l,l]\times \{-m\}$, 
 and let $z_1,z_2,\ldots, z_s$ be a black 8-path from  $\{1\}\times \{1\}$ to  $[-l,l]\times \{m\}$. Let 
 \begin{equation*}
 \begin{split}
 D=&\{w_1,w_2,\ldots,w_p, w_p,\ldots,w_2,w_1, x_1,x_2,\ldots,x_q, x_q, \ldots, x_2,x_1, y_1,y_2,\ldots, y_r, \\
 & y_r,\ldots, y_2,  y_1, z_1,z_2,\ldots, z_s,z_s,\ldots, z_2,z_1\}.
 \end{split}
 \end{equation*} 
 For $t\in D$, if $t=(a,b)$, we let $\psi_1(t)=a$ and $\psi_2(t)=b$. This works.
\end{proof}

In the rest of the section we prove Proposition \ref{path}.
Let $[r], [s]$, and $[t]$ be given. Let $\alpha\colon [s]\to [r]$ and $\beta\colon [t]\to [r]$ be relation preserving maps (not necessarily onto).
For the chessboard $ \mathbb{C}=[s]\times [t]$,  we let $(i,j)$ to be black if and only if $ \alpha(i)=\beta(j)$. 
In the lemma below we collect  amalgamation results we will use later.
\begin{lemma}\label{amalg}

\noindent a) If $\alpha(1)=1$, $\beta(1)=\beta(t)=1$, and $\mbox{rng}(\beta)\subseteq \mbox{rng}(\alpha)$, then there is a black 8-path from $(1,1)$ to
$(1,t)$. 

\noindent b) If  $\alpha(1)=1$, $\beta(1)=1$, and $\mbox{rng}(\beta)\subseteq\mbox{rng}(\alpha)$, then there is a black 8-path from  $(1,1)$ to $[s]\times\{t\}$.

\end{lemma}

We will write a careful proof of  part a). A proof of b) is very similar, and is left to the reader.

\begin{proof}
We show that there is a black  8-path from $(1,1)$ to $(1,t)$. 
 For this, via the Steinhaus' chessboard theorem, it is enough to 
show that there is no white 4-path from $\ov{(1,1)(1,t)}$ to $\ov{(1,t)(1,1)}$. Take a 4-path $(a_1,b_1),(a_2,b_2),\ldots,(a_n,b_n)$ from $\ov{(1,1)(1,t)}$ to 
$\ov{(1,t)(1,1)}$.
We show that for some $i$, $(a_i,b_i)$ is black, that is, $\alpha(a_i)=\beta(b_i)$. 
Define $h\colon \{1,2,\ldots, n\}\to \mathbb{Z}$, where $\mathbb{Z}$ is the set of integers,  by 
$h(i)=\alpha(i)-\beta(i)$. The function $h$ has an important `continuity' property:  for every $i$, $|h(i+1)-h(i)|\leq 1$. We want to find $i$ such that $h(i)=0$.
We will consider three cases: $(a_n,b_n)\in [s]\times \{1\}$, $(a_n,b_n)\in [s]\times \{t\}$, and $(a_n,b_n)\in \{s\}\times [t]$. 

First, let $(a_n,b_n)\in [s]\times \{1\}$.  Since $\alpha(a_1)=1$,  we have $h(1)\leq 0$, and since $\beta(b_n)=1$, we have $h(n)\geq 0$. Therefore, by the `continuity' property,
for some $i$, $h(i)=0$. In the case when $(a_n,b_n)\in [s]\times \{t\}$, for the same reason, there is $i$ such that $h(i)=0$.
Suppose now that $(a_n,b_n)\in \{s\}\times [t]$.  Let $[r_0]=\mbox{rng}(\alpha)$.  Let $x,y\in [s]$ be such that $\alpha(x)=1$ and $\alpha(y)=r_0$. Take $(a_{i_0}, b_{i_0})$
 such that $a_{i_0}=x$ and take $(a_{j_0}, b_{j_0})$ such that $a_{j_0}=y$.
 Since  $\mbox{rng}(\beta)\subseteq \mbox{rng}(\alpha)$, we have $\beta(b_{j_0})\leq r_0$, and therefore $h(j_0)\geq 0$. Since also 
 $h(i_0)\leq 0$, for some $i$ we have $h(i)=0$.

\end{proof}

\begin{remark}\label{remsla}
The Steinhaus' chessboard theorem was used by  Solecki to prove the AP of the family $\g$ (unpublished). His proof is much simpler than the one presented in \cite{IS}.
We give here a sketch of Solecki's proof (with his permission). Let $A=[k],B=[l],C=[m]\in\g$ and epimorphisms $\alpha\colon B\to A$ and $\beta\colon B\to A$ be given.  
We want to find $D\in\g$, $\gamma \colon D\to B$, and $\delta\colon D\to C$ such that $\alpha\circ\gamma=\beta\circ\delta$.
Consider the chessboard $[l]\times [m]$. Let $(i,j)$ be black if and only if $\alpha(i)=\beta(j)$. An argument similar to the one we used in the proof of Lemma \ref{amalg}
shows that there is no white 4-path from $\{1\}\times [m]$ to $\{l\}\times [m]$, and there is no white 4-path from $[l]\times\{1\}$ to  $[l]\times\{m\}$.
Therefore, by the  Steinhaus' chessboard theorem, there are a black 8-path from  $[l]\times\{1\}$ to  $[l]\times\{m\}$ and  a black 8-path from 
$\{1\}\times [m]$ to $\{l\}\times [m]$.
We can combine these two 8-paths into one 8-path $x_1,x_2,\ldots, x_n$
with the following property:
for every $a\in [l]$ there is $b\in [m]$ such that for some $i$, $(a,b)=x_i$, and
for every $b\in [m]$ there is $a\in [l]$  such that for some $i$, $(a,b)=x_i$.
Let $D=[n]$. For $i\in [n]$, if $x_i=(a,b)$, we let $\gamma(i)=a$ and $\delta(i)=b$. This works.
\end{remark}

Let $A,B,C$ and $\phi_1,\phi_2$ be as in the hypotheses of Proposition \ref{path}. We pick points in $B=[-l,-1]\cup[1,l]$:
\[
s_{-p}<s'_{-p+1}<\ldots<s'_{-2}<s_{-2}<s'_{-1}<s_{-1}<s'_0<s_0<s'_1<s_1<s'_2<s_2<\ldots<s'_p,
\]
for some $p$, so that $s_0=1$, $s'_0=-1$, $s'_p=l$, $s_{-p}=-l$, for every $-p<i<p$ we have $s_i=s'_i +1$,  for every $-p\leq i<p$ the epimorphism $\phi_1$ 
assumes only positive values or assumes only negative values in the interval $[s_i,s'_{i+1}]$,
 and for every $-p<i<p$, $\phi_1(s'_i)$ and $\phi_1(s_i)$ have opposite signs.
 
 Notice that since $\phi_1(j)=-\phi_1(-j)$, $j\in B$, we have  for every $-p\leq i<p$, $s_i=-s'_{-i}$. Notice also that if $\phi_1(s_0)>0$ then for any even number $-p\leq i<p$  
  we have $\phi_1(s_i)=\phi_1(s'_{i+1})=1$ and for any odd number  $-p\leq i<p$ we have $\phi_1(s_i)=\phi_1(s'_{i+1})=-1$.
 On the other hand, if $\phi_1(s_0)<0$ then for any even number $-p\leq i<p$  we have $\phi_1(s_i)=\phi_1(s'_{i+1})=-1$ and for any odd number 
  $-p\leq i<p$  we have $\phi_1(s_i)=\phi_1(s'_{i+1})=1$.

 Similarly, we pick points in $C=[-m,-1]\cup[1,m]$:
\[
t_{-q}<t'_{-q+1}<\ldots<t'_{-2}<t_{-2}<t'_{-1}<t_{-1}<t'_0<t_0<t'_1<t_1<t'_2<t_2<\ldots<t'_q,
\]
for some $q$, so that $t_0=1$, $t'_0=-1$, $t'_q=m$, $t_{-q}=-m$, for every $-q<i<q$ we have $t_i=t'_i +1$,  for every $-q\leq i<q$ the epimorphism $\phi_2$ 
assumes only positive values or assumes only negative values in the interval $[t_i,t'_{i+1}]$,
 and for every $-q<i<q$, $\phi_2(t'_i)$ and $\phi_2(t_i)$ have opposite signs.
 
 Notice that since $\phi_2(j)=-\phi_2(-j)$, $j\in C$, we have  for every $-q\leq i<q$, $t_i=-t'_{-i}$. Notice also that if $\phi_2(t_0)>0$ then for any even number $-q\leq i<q$
  we have $\phi_2(t_i)=\phi_2(t'_{i+1})=1$ and for any odd number $-q\leq i<q$  we have $\phi_2(t_i)=\phi_2(t'_{i+1})=-1$.
 On the other hand, if $\phi_2(t_0)<0$ then for any even number $-q\leq i<q$
  we have $\phi_2(t_i)=\phi_2(t'_{i+1})=-1$ and for any odd number $-q\leq i<q$  we have $\phi_2(t_i)=\phi_2(t'_{i+1})=1$.

\bigskip

We define a graph $\mathcal{G}_0$. Let the set of vertices in  $\mathcal{G}_0$ be equal to the set $V=[-p,p]\times [-q,q]$. For  $(a,b),(c,d)\in V$, if they are not 4-adjacent,
they will not be connected by an edge. For every $-p\leq i<p$ and $-q\leq j<q$ such that $\phi_1$ has the same sign on $[s_i,s'_{i+1}]$ as $\phi_2$ has 
on $[t_j, t'_{j+1}]$, if $\mbox{rng}(\phi_1\restriction [s_i,s'_{i+1}])\subseteq \mbox{rng}(\phi_2\restriction [t_j, t'_{j+1}])$, we put an edge between $(i,j) $ and $(i+1,j) $ and we 
put an edge between $(i,j+1)$ and $(i+1,j+1)$.
 If $\mbox{rng}(\phi_2\restriction [t_j, t'_{j+1}])\subseteq \mbox{rng}(\phi_1\restriction [s_i,s'_{i+1}])$, 
we put an edge between $(i,j)$ and $(i,j+1)$ and we put an edge between $(i+1, j)$ and $(i+1,j+1)$.
 We have just defined $\mathcal{G}_0$.

 Let $x_1,x_2,\ldots, x_n$ be a path in the chessboard $[-p,p]\times [-q,q]$. We say that it is {\em interior} if  $x_1,\ldots, x_{n-1}\notin\mbox{Bd}([-p,p]\times [-q,q])$.
 \begin{lemma}\label{long}
 \begin{enumerate}
 \item The existence of an interior path  in $\mathcal{G}_0$  (in the graph-theoretic sense) from $(0,0)$ to $[-p,p]\times \{-q\}$ implies 
 the existence of a black 8-path    in $B\times C$ from  $\{-1,1\}\times \{-1,1\}$ to $[-l,l]\times \{-m\}$.
 \item The existence of an interior path in $\mathcal{G}_0$    from $(0,0)$ to $[-p,p]\times \{q\}$ implies 
 the existence of a black 8-path in $B\times C$ from  $\{-1,1\}\times \{-1,1\}$ to  $[-l,l]\times \{m\}$.
 \item The existence of an interior path  in the graph $\mathcal{G}_0$ 
 from $(0,0)$ to $\{-p\}\times [-q,q]$, implies the existence of  a black 8-path  in $B\times C$ from $\{-1,1\}\times \{-1,1\}$ to $\{-l\}\times [-m,m]$. 
 \item The existence of an interior path  in $\mathcal{G}_0$   from $(0,0)$ to $\{p\}\times [-q,q]$ implies 
 the existence of a black 8-path  in $B\times C$ from  $\{-1,1\}\times \{-1,1\}$ to $\{l\}\times [-m,m]$. 
 \end{enumerate}
 \end{lemma}

\begin{proof}
This follows from Lemma \ref{amalg} and from the definition of $\mathcal{G}_0$.
\end{proof}
 
To finish the proof of Proposition \ref{path}, we have to show that there are paths in $\mathcal{G}_0$, as in the hypotheses of Lemma \ref{long}.
 
 Define $\mathcal{G}_1$ to be a subgraph of $\mathcal{G}_0$ such that for every $-p\leq i<p$ and $-q\leq j<q$
 such that $\phi_1$ has the same sign on $[s_i,s'_{i+1}]$ as $\phi_2$ has on $[t_j, t'_{j+1}]$
 if $\mbox{rng}(\phi_2\restriction [t_j, t'_{j+1}])= \mbox{rng}(\phi_1\restriction [s_i,s'_{i+1}])$, we  delete the edge between $(i,j) $ and $(i+1,j) $, we 
 delete the edge between $(i,j+1)$ and $(i+1,j+1)$, keep the edge between $(i,j)$ and $(i,j+1)$, and we keep the edge between $(i+1, j)$ and $(i+1,j+1)$.

Define $\mathcal{G}_2$  to be a subgraph of $\mathcal{G}_0$ such that for every $-p\leq i<p$ and $-q\leq j<q$
such that $\phi_1$ has the same sign on $[s_i,s'_{i+1}]$ as $\phi_2$ has on $[t_j, t'_{j+1}]$
 if $\mbox{rng}(\phi_2\restriction [t_j, t'_{j+1}])= \mbox{rng}(\phi_1\restriction [s_i,s'_{i+1}])$, we  keep the edge between $(i,j) $ and $(i+1,j) $, we 
 keep the edge between $(i,j+1)$ and $(i+1,j+1)$, delete the edge between $(i,j)$ and $(i,j+1)$, and we delete the edge between $(i+1, j)$ and $(i+1,j+1)$.
 
 \begin{remark}\label{sym}
In the graph  $\mathcal{G}_0$, $(a,b)$ and $(c,d)$ are connected by an edge if and only if $(-a,-b)$ and $(-c,-d)$ are connected by an edge.
 The same conclusion holds for graphs $\mathcal{G}_1$ and $\mathcal{G}_2$.
 \end{remark}
 
 The existence of the required paths in $\mathcal{G}_0$ will follow from the lemma below.
 \begin{lemma}\label{graf}
 \begin{enumerate}
 \item  In the graph $\mathcal{G}_1$ there is an interior path   from $(0,0)$ to $[-p,p]\times \{-q\}$
 and there is an interior path    from $(0,0)$ to $[-p,p]\times \{q\}$.
 \item  In the graph $\mathcal{G}_2$ there is an interior path  from $(0,0)$ to $\{-p\}\times [-q,q]$ 
 and  there is an interior path from $(0,0)$ to $\{p\}\times [-q,q]$. 
 \end{enumerate}
 \end{lemma}
 
 \begin{proof}
 We show (1). The proof of (2) is similar. 
 
 \noindent {\bf{Claim 1.}} For every $-p<i<p$, $-q<j<q$  there are exactly two edges that end in $(i,j)$.
 \begin{proof}
Fix $(i,j)$. Either  $\phi_1$ has the same sign on $[s_i,s'_{i+1}]$ as $\phi_2$ has on $[t_j,t'_{j+1}]$ and  $\phi_1$ has the same sign on $[s_{i-1},s'_{i}]$  as $\phi_2$ has on $[t_{j-1},t'_{j}]$, or $\phi_1$ has the same sign on $[s_i,s'_{i+1}]$ as $\phi_2$ has on $[t_{j-1},t'_{j}]$ and  $\phi_1$ has the same sign on $[s_{i-1},s'_{i}]$ as $\phi_2$ has on $[t_{j},t'_{j+1}]$ (exactly one of these two possibilities hold).
 Say, the former holds. Since $\phi_1$ has the same sign on $[s_i,s'_{i+1}]$ as $\phi_2$ has on $[t_j,t'_{j+1}]$, by the definition of $\mathcal{G}_1$, there is an edge between 
 $(i,j)$ and $(i+1,j)$, or there is an edge between $(i,j)$ and $(i,j+1)$ (exactly one of these two possibilities hold). Since $\phi_1$ has the same sign on $[s_{i-1},s'_{i}]$ as $\phi_2$ has on $[t_{j-1},t'_{j}]$,
 there is an edge between $(i,j)$ and $(i-1,j)$, or there is an edge between $(i,j)$ and $(i,j-1)$ (exactly one of these two possibilities hold). This finishes the proof of the claim.
\end{proof}
 
  \noindent {\bf{Claim 2.}} There are no loops passing through $(0,0)$ in the graph $\mathcal{G}_1$.
  \begin{proof}
Suppose towards a contradiction that $(0,0)=(a_0,b_0),(a_1,b_1),\ldots, (a_n,b_n)=(0,0)$,  where $(a_0,b_0),(a_1,b_1),\ldots, (a_{n-1},b_{n-1})$ are pairwise different,
is a loop. By Claim~1 and Remark \ref{sym}, $(a_{n-1},b_{n-1})=(-a_1,-b_1)$, $(a_{n-2},b_{n-2})=(-a_2,-b_2)$,...
Hence, if $n$ is even, we have $(a_{\frac{n}{2}},b_{\frac{n}{2}}) =(-a_{\frac{n}{2}},-b_{\frac{n}{2}}) $. This implies $(a_{\frac{n}{2}},b_{\frac{n}{2}}) =(0,0)$ and 
contradicts the assumption that $ (a_{\frac{n}{2}},b_{\frac{n}{2}}) \neq (a_0,b_0)$.
If $n$ is odd, we get $(a_{\frac{n-1}{2}},b_{\frac{n-1}{2}})= (-a_{\frac{n+1}{2}},-b_{\frac{n+1}{2}}) $. Since $(a_{\frac{n-1}{2}},b_{\frac{n-1}{2}})$
and $ (a_{\frac{n+1}{2}},b_{\frac{n+1}{2}}) $ are adjacent, and since $(a_{\frac{n-1}{2}},b_{\frac{n-1}{2}})\neq (0,0) $, we again arrive at a contradiction.
  \end{proof}
  
  Notice that Claim 1 and Claim 2 already imply that there is a path from $(0,0)$ to the boundary of $[-p,p]\times [-q,q]$.
  
  Let $-p\leq i_0<p$ be such that $\mbox{rng}(\phi_1\restriction [s_{i_0}, s'_{i_0+1}])=[1,k]$ (recall that $A=[-k,-1]\cup [1,k]$). Then  
  $\mbox{rng}(\phi_1\restriction [s_{-(i_0+1)}, s'_{-i_0}])=[-k,-1]$. Therefore, 
  for every $-q\leq j<q$,
if $\phi_1$ has the same sign on $[s_{i_0},s'_{i_0+1}]$ as $\phi_2$ has on $[t_j, t'_{j+1}]$, then
 $\mbox{rng}(\phi_2\restriction [t_j, t'_{j+1}])\subseteq \mbox{rng}(\phi_1\restriction [s_{i_0},s'_{i_0+1}])$, and if
 $\phi_1$ has the same sign on $ [s_{-(i_0+1)}, s'_{-i_0}]$ as $\phi_2$ has on $[t_j, t'_{j+1}]$, then 
  $\mbox{rng}(\phi_2\restriction [t_j, t'_{j+1}])\subseteq \mbox{rng}(\phi_1\restriction [s_{-(i_0+1)}, s'_{-i_0}]$.
  Therefore, for every $j$, there is no edge between $(i_0,j)$ and $(i_0+1,j)$ and there is no edge between $(-(i_0+1),j)$ and $(-i_0,j)$. 
  
  This implies that the path from $(0,0)$ to the boundary of $[-p,p]\times [-q,q]$ is in fact from $(0,0)$ to $[-i_0,i_0]\times \{-q,q\}$.
   Since   $\phi_1(i)=-\phi_1(-i)$ and  $\phi_2(i)=-\phi_2(-i)$, the existence of  a path  from $(0,0)$ to $[-p,p]\times \{-q\}$ 
 is equivalent to the existence of a path  from $(0,0)$ to $[-p,p]\times \{q\}$. Therefore, there exist   paths  from $(0,0)$ to $[-p,p]\times \{-q\}$ 
and  from $(0,0)$ to $[-p,p]\times \{q\}$. 
Without loss of generality, only the last element of each path is in  $\mbox{Bd}([-p,p]\times [-q,q])$.
 
 \end{proof}
 
\begin{example}\label{exa}
We illustrate the proof of Theorem \ref{maint} on an example. Let $A=[-3,-1]\cup [1,3], B=[-8,-1]\cup [1,8], C=[-9,-1]\cup [1,9]$.
Let $\phi_1\colon B\to A$ be given by $\phi_1(1)=1$, $\phi_1(2)=2$, $\phi_1(3)=1$, $\phi_1(4)=-1$, $\phi_1(5)=1$, $\phi_1(6)=1$, $\phi_1(7)=2$, $\phi_1(8)=3$,
and $\phi_1(-i)=\phi_1(i)$ for $i\in[1,8]$.
Let $\phi_2\colon C\to A$ be given by $\phi_2(1)=-1$, $\phi_2(2)=-2$, $\phi_2(3)=-1$, $\phi_2(4)=-2$, $\phi_2(5)=-3$, $\phi_2(6)=-2$, $\phi_2(7)=-1$, $\phi_2(8)=1$,
$\phi_2(9)=2$,
and $\phi_2(-i)=\phi_2(i)$ for $i\in[1,9]$.

We consider the chessboard $B\times C$, where $(i,j)$ is black if and only if $\phi_1(i)=\phi_2(j)$ (Figure 1).

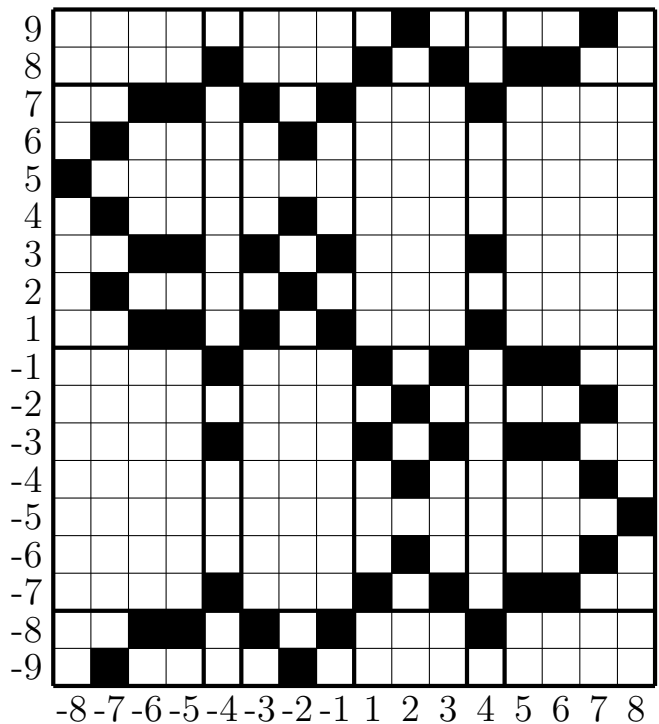
\begin{figure}
\begin{tikzpicture}[scale=0.5]
\draw [ultra thick] (-8,0)--(8,0);                            
\draw (-8,1)--(8,1);
\draw (-8,2)--(8,2);
\draw (-8,3)--(8,3);
\draw (-8,4)--(8,4);
\draw (-8,5)--(8,5);
\draw (-8,6)--(8,6);
\draw  [ultra thick] (-8,7)--(8,7);
\draw (-8,8)--(8,8);
\draw  [ultra thick] (-8,9)--(8,9);
\draw (-8,-1)--(8,-1);
\draw (-8,-2)--(8,-2);
\draw (-8,-3)--(8,-3);
\draw (-8,-4)--(8,-4);
\draw (-8,-5)--(8,-5);
\draw (-8,-6)--(8,-6);
\draw  [ultra thick] (-8,-7)--(8,-7);
\draw (-8,-8)--(8,-8);
\draw  [ultra thick] (-8,-9)--(8,-9);

\draw  [ultra thick] (0,-9)--(0,9);
\draw (1,-9)--(1,9);
\draw (2,-9)--(2,9);
\draw  [ultra thick] (3,-9)--(3,9);
\draw  [ultra thick] (4,-9)--(4,9);
\draw (5,-9)--(5,9);
\draw (6,-9)--(6,9);
\draw (7,-9)--(7,9);
\draw  [ultra thick] (8,-9)--(8,9);
\draw (-1,-9)--(-1,9);
\draw (-2,-9)--(-2,9);
\draw  [ultra thick] (-3,-9)--(-3,9);
\draw  [ultra thick] (-4,-9)--(-4,9);
\draw (-5,-9)--(-5,9);
\draw (-6,-9)--(-6,9);
\draw (-7,-9)--(-7,9);
\draw  [ultra thick] (-8,-9)--(-8,9);

\draw [fill] (0,7) rectangle (1,8);  \draw [fill] (0,-7) rectangle (-1,-8);
\draw [fill] (1,8) rectangle (2,9);  \draw [fill] (-1,-8) rectangle (-2,-9);
\draw [fill] (2,7) rectangle (3,8);  \draw [fill] (-2,-7) rectangle (-3,-8);
\draw [fill] (3,0) rectangle (4,1);  \draw [fill] (-3,0) rectangle (-4,-1);
\draw [fill] (3,2) rectangle (4,3);  \draw [fill] (-3,-2) rectangle (-4,-3);
\draw [fill] (3,6) rectangle (4,7);  \draw [fill] (-3,-6) rectangle (-4,-7);
\draw [fill] (4,7) rectangle (5,8);  \draw [fill] (-4,-7) rectangle (-5,-8);
\draw [fill] (5,7) rectangle (6,8);  \draw [fill] (-5,-7) rectangle (-6,-8);
\draw [fill] (6,8) rectangle (7,9);  \draw [fill] (-6,-8) rectangle (-7,-9);

\draw [fill] (0,-1) rectangle (1,0);  \draw [fill] (0,1) rectangle (-1,0);
\draw [fill] (0,-3) rectangle (1,-2);  \draw [fill] (0,3) rectangle (-1,2);
\draw [fill] (0,-7) rectangle (1,-6);  \draw [fill] (0,7) rectangle (-1,6);
\draw [fill] (2,-1) rectangle (3,0);  \draw [fill] (-2,1) rectangle (-3,0);
\draw [fill] (2,-3) rectangle (3,-2);  \draw [fill] (-2,3) rectangle (-3,2);
\draw [fill] (2,-7) rectangle (3,-6);  \draw [fill] (-2,7) rectangle (-3,6);
\draw [fill] (4,-1) rectangle (5,0);  \draw [fill] (-4,1) rectangle (-5,0);
\draw [fill] (4,-3) rectangle (5,-2);  \draw [fill] (-4,3) rectangle (-5,2);
\draw [fill] (4,-7) rectangle (5,-6);  \draw [fill] (-4,7) rectangle (-5,6);
\draw [fill] (5,-1) rectangle (6,0);  \draw [fill] (-5,1) rectangle (-6,0);
\draw [fill] (5,-3) rectangle (6,-2);  \draw [fill] (-5,3) rectangle (-6,2);
\draw [fill] (5,-7) rectangle (6,-6);  \draw [fill] (-5,7) rectangle (-6,6);
\draw [fill] (0,-1) rectangle (1,0);  \draw [fill] (0,1) rectangle (-1,0);
\draw [fill] (0,-3) rectangle (1,-2);  \draw [fill] (0,3) rectangle (-1,2);
\draw [fill] (3,-8) rectangle (4,-7);  \draw [fill] (-3,8) rectangle (-4,7);
\draw [fill] (7,-5) rectangle (8,-4);  \draw [fill] (-7,5) rectangle (-8,4);
\draw [fill] (1,-6) rectangle (2,-5);  \draw [fill] (-1,6) rectangle (-2,5);
\draw [fill] (1,-4) rectangle (2,-3);  \draw [fill] (-1,4) rectangle (-2,3);
\draw [fill] (1,-2) rectangle (2,-1);  \draw [fill] (-1,2) rectangle (-2,1);
\draw [fill] (6,-6) rectangle (7,-5);  \draw [fill] (-6,6) rectangle (-7,5);
\draw [fill] (6,-4) rectangle (7,-3);  \draw [fill] (-6,4) rectangle (-7,3);
\draw [fill] (6,-2) rectangle (7,-1);  \draw [fill] (-6,2) rectangle (-7,1);

\node[below] at (0.5,-9) {\large{1}};   \node[below] at (-0.5,-9) {\large{-1}};
\node[below] at (1.5,-9) {\large{2}};   \node[below] at (-1.5,-9) {\large{-2}};
\node[below] at (2.5,-9) {\large{3}};   \node[below] at (-2.5,-9) {\large{-3}};
\node[below] at (3.5,-9) {\large{4}};   \node[below] at (-3.5,-9) {\large{-4}};
\node[below] at (4.5,-9) {\large{5}};   \node[below] at (-4.5,-9) {\large{-5}};
\node[below] at (5.5,-9) {\large{6}};   \node[below] at (-5.5,-9) {\large{-6}};
\node[below] at (6.5,-9) {\large{7}};   \node[below] at (-6.5,-9) {\large{-7}};
\node[below] at (7.5,-9) {\large{8}};   \node[below] at (-7.5,-9) {\large{-8}};

\node[left] at (-8,0.5) {\large{1}};   \node[left] at (-8,-0.5) {\large{-1}};
\node[left] at (-8,1.5) {\large{2}};   \node[left] at (-8,-1.5) {\large{-2}};
\node[left] at (-8,2.5) {\large{3}};   \node[left] at (-8,-2.5) {\large{-3}};
\node[left] at (-8,3.5) {\large{4}};   \node[left] at (-8,-3.5) {\large{-4}};
\node[left] at (-8,4.5) {\large{5}};   \node[left] at (-8,-4.5) {\large{-5}};
\node[left] at (-8,5.5) {\large{6}};   \node[left] at (-8,-5.5) {\large{-6}};
\node[left] at (-8,6.5) {\large{7}};   \node[left] at (-8,-6.5) {\large{-7}};
\node[left] at (-8,7.5) {\large{8}};   \node[left] at (-8,-7.5) {\large{-8}};
\node[left] at (-8,8.5) {\large{9}};   \node[left] at (-8,-8.5) {\large{-9}};

\end{tikzpicture}
\caption{The chessboard $B\times C$}
\end{figure}

Therefore, $s_{-3}=-8$, $s'_{-2}=-5$, $s_{-2}=-4$, $s'_{-1}=-4$, $s_{-1}=-3$, $s'_0=-1$
$s_0=1$, $s'_1=3$, $s_1=4$, $s'_2=4$, $s_2=5$, $s'_3=8$, and
$t_{-2}=-9$, $t'_{-1}=-8$, $t_{-1}=-7$, $t'_0=-1$,
$t_0=1$, $t'_1=7$, $t_1=8$, $t'_2=9$.

Moreover, $\mbox{rng}(\phi_1\restriction [s_{-3},s'_{-2}])=[-3,-1]$, $\mbox{rng}(\phi_1\restriction [s_{-2},s'_{-1}])=[1]$, $\mbox{rng}(\phi_1\restriction [s_{-1},s'_{0}])=[-2,-1]$,
$\mbox{rng}(\phi_1\restriction [s_0,s'_{1}])=[1,2]$, $\mbox{rng}(\phi_1\restriction [s_1,s'_{2}])=[-1]$, $\mbox{rng}(\phi_1\restriction [s_2,s'_{3}])=[1,3]$, and
$\mbox{rng}(\phi_2\restriction [t_{-2},t'_{-1}])=[-2,-1]$,  $\mbox{rng}(\phi_2\restriction [t_{-1},t'_{0}])=[1,3]$,
  $\mbox{rng}(\phi_2\restriction [t_{0},t'_{1}])=[-3,-1]$,  $\mbox{rng}(\phi_2\restriction [t_{1},t'_{2}])=[1,2]$.
  
  Hence, graphs $\mathcal{G}_1$ and $\mathcal{G}_2$ are as in Figure 2.

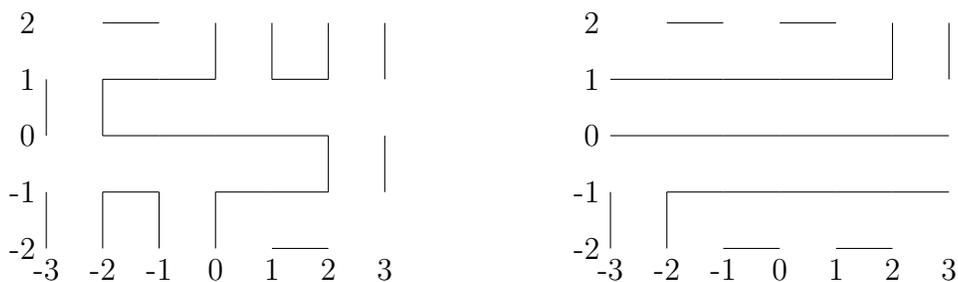
\begin{figure}
\begin{tikzpicture}[scale=0.75]
\draw (-3,-2) -- (-3,-1);
\draw (-2,-2) -- (-2,-1);
\draw (-1,-2) -- (-1,-1);
\draw (0,-2) -- (0,-1);
\draw (2,-1) -- (2,0);
\draw (3,-1) -- (3,0);
\draw (3,2) -- (3,1);
\draw (2,2) -- (2,1);
\draw (1,2) -- (1,1);
\draw (0,2) -- (0,1);
 \draw (-2,1) -- (-2,0);
 \draw (-3,1) -- (-3,0);

 \draw (1,-2) -- (2,-2);
 \draw (1,-1) -- (2,-1);
 \draw (-2,-1) -- (-1,-1);
 \draw (-2,0) -- (-1,0);
 \draw (0,-1) -- (1,-1);
 \draw (0,0) -- (1,0);
 \draw (-1,2) -- (-2,2);
 \draw (-1,1) -- (-2,1);
 \draw (2,1) -- (1,1);
 \draw (2,0) -- (1,0);
 \draw (0,1) -- (-1,1);
 \draw (0,0) -- (-1,0);


 \node[below] at (0,-2) {0};   
 \node[below] at (1,-2) {1};   \node[below] at (-1,-2) {-1};
 \node[below] at (2,-2) {2};   \node[below] at (-2,-2) {-2};
 \node[below] at (3,-2) {3};   \node[below] at (-3,-2) {-3};

 \node[left] at (-3,0)  {0};  
 \node[left] at (-3,1)  {1}; 
 \node[left] at (-3,2)  {2}; 
 \node[left] at (-3,-1)  {-1}; 
 \node[left] at (-3,-2)  {-2}; 

\draw (13,1) -- (13,2);
\draw (7,-2) -- (7,-1);
\draw (8,-2) -- (8,-1);
\draw(11,1) -- (12,1);
 \draw (9,-2) -- (10,-2);
\draw (9,-1) -- (10,-1);
\draw (11,-2) -- (12,-2);
\draw (11,-1) -- (12,-1);
\draw (8,-1) -- (9,-1);
\draw (8,0) -- (9,0);
\draw (10,-1) -- (11,-1);
\draw (10,0) -- (11,0);
\draw (12,-1) -- (13,-1);
\draw (12,0) -- (13,0);
\draw (11,2) -- (10,2);
\draw (11,1) -- (10,1);
\draw (9,2) -- (8,2);
\draw (9,1) -- (8,1);
\draw (12,1) -- (12,2);
\draw (12,0) -- (11,0);
\draw (10,1) -- (9,1);
\draw (10,0) -- (9,0);
\draw (8,1) -- (7,1);
\draw (8,0) -- (7,0);


\node[below] at (10,-2) {0};   
\node[below] at (11,-2) {1};   \node[below] at (9,-2) {-1};
\node[below] at (12,-2) {2};   \node[below] at (8,-2) {-2};
\node[below] at (13,-2) {3};   \node[below] at (7,-2) {-3};

\node[left] at (7,0)  {0};  
\node[left] at (7,1)  {1}; 
\node[left] at (7,2)  {2}; 
\node[left] at (7,-1)  {-1}; 
\node[left] at (7,-2)  {-2}; 
\end{tikzpicture}
\caption{Graphs $\mathcal{G}_1$ (left) and $\mathcal{G}_2$ (right)} 
\end{figure}

\end{example}

\appendix
\section{}
The purpose of this appendix is to present a criterium for the automorphism group of a projective Fra\"{i}ss\'{e} limit acting  by conjugacy to have a dense orbit, and 
to present a criterium for the automorphism group of a projective Fra\"{i}ss\'{e} limit  acting  by conjugacy to have a comeager orbit. These criteria and their proofs are  analogs of Theorems 2.1
 and 3.4 given by Kechris and Rosendal \cite{KR} in the context of  (injective)  Fra\"{i}ss\'{e} limits. However, we point out that we will work with surjective relations rather 
 than with partial functions, and our criteria are analogs but not  dualizations of the corresponding criteria in \cite{KR}. It seems that working with surjective relations rather 
 than with partial functions  makes calculations simpler in the context of  projective  Fra\"{i}ss\'{e} limits.
 We hope that  many new interesting projective Fra\"{i}ss\'{e} limits will be discovered, and these criteria will be useful for them.

Let $\g$ be a countable projective Fra\"{i}ss\'{e} family in a language $L_0$. Let $\p$ be the projective Fra\"{i}ss\'{e} limit of $\g$. Let $\aut$ be the automorphism group of $\p$ equipped with the topology inherited from the homeomorphism group $H(\p)$. The $\aut$ acts on itself by conjugacy  and let $H$ be a closed invariant subspace of $\aut$.
Define \begin{equation*}
\begin{split}
\f= & \{ (A, s^A)\colon A\in\g \mbox{ and } \exists\phi\colon\p\to A \exists f\in H \mbox{ such that}\\
 &\phi\colon (\p,f)\to(A,s^A) \mbox{ is an epimorphism}\}.
\end{split}
\end{equation*}

\begin{theorem}
Let $\g$, $\p$, $H$, and $\f$ be as above.
 Then $\f$ has the JPP if and only if $\aut\acts H$ has a dense orbit.
\end{theorem}

\begin{proof}
The proof that the JPP implies that $\aut\acts H$ has a dense orbit is the same as in the special case (Proposition \ref{dense}).

We show the converse.
Take $(A,s^A), (B,s^B)\in\f$. We find $(C,s^C)\in\f$ such that there are epimorphisms from $(C,s^C)$ onto $(A,s^A)$ and onto $(B,s^B)$. Take any epimorphism
$\phi\colon\p\to A$. Take $f\in [\phi,s^A]$ that has a dense orbit. Take any epimorphism $\psi\colon\p\to B$.
Let $g\in\aut$ be such that $gfg^{-1}\in [\psi, s^B]$.  Let $C$ be a partition of $\p$ that refines partitions $\phi^{-1}(A)$ and $g(\psi^{-1}(B))$, 
and moreover $r^\p$ restricted to $C$ is in $\g$. (To achieve this last requirement on $C$, we use (L2).) Let $\hat{\phi}$ be the natural projection from $\p$
to $C$. We let $s^C(c,d)$ if and only if there are $x,y\in\p$ such that $\hat{\phi}(x)=c$,   $\hat{\phi}(y)=d$, and $f(x)=y$. Clearly, the natural projection $\bar{\phi}$ from
$(C,s^C)$ onto $(A,s^A)$ is an epimorphism. Let $\hat{\psi}=\hat{\phi}\circ g^{-1}$. Let $\bar{\psi}$ be the natural projection from 
$(g^{-1}(C), g^{-1}(s^C))$ to $(B,s^B)$. Since there are $x,y\in\p$ such that $\hat{\phi}(x)=c$,   $\hat{\phi}(y)=d$, and $f(x)=y$ if and only if
there are $x,y\in\p$ such that $\hat{\psi}(x)=g^{-1}(c)$,   $\hat{\psi}(y)=g^{-1}(d)$, and $gfg^{-1}(x)=y$, this projection is an epimorphism.

\end{proof}

We say that a family $\f$ of topological $L$-structures has  {\em the weak amalgamation property}, or {\em the WAP}, if for every $A\in\f$
there is $B\in\f$ and an epimorphism $\phi\colon B\to A$ such that
  for any $C_1,C_2\in\f$ and any epimorphisms $\phi_1\colon C_1\to B$ and $\phi_2\colon C_2\to B$, there exist $D\in\f$,
 $\phi_3\colon D\to C_1$, and $\phi_4\colon D\to C_2$ such that $\phi\circ\phi_1\circ \phi_3=\phi\circ\phi_2\circ \phi_4$.

\begin{theorem}
Let $\g$, $\p$, $H$, and $\f$ be as above.
 Then $\f$ has the JPP and the WAP if and only if $\aut\acts H$ has a comeager orbit.
\end{theorem}

\begin{proof}
The proof that the JPP and the CAP imply that $\aut\acts H$ has a comeager orbit is the same as in the special case (Proposition \ref{comeager}).
To show that the JPP and the WAP imply that $\aut\acts H$ has a comeager orbit we have to make small modifications.
Namely, since we have only WAP and not necessarily CAP, we take the following definition of $E(\phi,s^A)$: 
\begin{equation*}
\begin{split}
 E(\phi,s^A)= & \{f\in H\colon \mbox{ if }    f\in[\phi,s^A] \mbox{ then for some } ((B,s^B),\psi,\bar{\psi}), \mbox{ an extension  }\\
  & \mbox{of } ((A,s^A),\phi) \mbox{ and }  ((B_m,s^{B_m}),\psi_m,\bar{\psi}_m) 
  \mbox{ an extension of  } ((B,s^B),\psi,\bar{\psi})\\
 & \mbox{ witnessing WAP, we have } f\in[\psi_m, s^{B_m}] \}. \\
\end{split}
 \end{equation*}

We now show the converse.
Let $(M,s^M)\in\f$. Let $\phi\colon\p\to M$ be an epimorphism.
Take $f\in [M,s^M]$ that has a  comeager orbit. By Proposition \ref{KR}, there is an open neighborhood $U$ of $f$ such that
$c(\phi,f)=\{gfg^{-1}\colon g\in [\phi,\mbox{id}_M]_{\aut}\}$ is dense in $U$. Therefore, there is $((N,s^N),\psi,\bar{\psi})$, an extension of $(M,s^M),\phi)$,
 such that  $c(\phi,f)$ is dense in $[\psi,s^N]$. We show that $(N,s^N)$ witnesses the WAP for $(M,s^M)$.
For this, take $(A,s^A), (B,s^B)\in\f$ and epimorphisms $\alpha\colon  (A,s^A)\to (N,s^N)$ and $\beta\colon (B,s^B)\to (N,s^N)$.
We find $(C,s^C)\in\f$ and epimorphisms $\gamma\colon(C,s^C)\to (A,s^A)$ and $\delta\colon(C,s^C)\to(B,s^B)$ such that
 $\bar{\psi}\circ\alpha\circ\gamma=\bar{\psi}\circ\beta\circ\delta$.
 
  Take an epimorphism
$\phi_1\colon\p\to A$ such that $\phi=\bar{\psi}\circ\alpha\circ\phi_1$. Take $f'\in [\phi_1,s^A]$, $f'=g_0fg_0^{-1}$ for some $
 g_0\in [\phi,\mbox{id}_M]_{\aut}$. Take an epimorphism $\phi_2\colon\p\to B$ such that $\phi=\bar{\psi}\circ\beta\circ\phi_2$.
Let $g\in [\phi,\mbox{id}_M]_{\aut}$ be such that $gf'g^{-1}\in [\phi_2, s^B]$.  Let $C$ be a partition of $\p$ that refines partitions $\phi_1^{-1}(A)$ and $g(\phi_2^{-1}(B))$, 
and moreover $r^\p$ restricted to $C$ is in $\g$. (To achieve this last requirement on $C$, we use (L2).) Let $\hat{\phi}_1$ be the natural projection from $\p$
to $C$. We let $s^C(c,d)$ if and only if there are $x,y\in\p$ such that $\hat{\phi}_1(x)=c$,   $\hat{\phi}_1(y)=d$, and $f'(x)=y$. Clearly, the natural projection $\bar{\phi}_1$ from
$(C,s^C)$ onto $(A,s^A)$ is an epimorphism. Let $\hat{\phi}_2=\hat{\phi}_1\circ g^{-1}$. Let $\bar{\phi}_2$ be the natural projection from 
$(g^{-1}(C), g^{-1}(s^C))$ to $(B,s^B)$. Since there are $x,y\in\p$ such that $\hat{\phi}_1(x)=c$,   $\hat{\phi}_1(y)=d$, and $f'(x)=y$ if and only if
there are $x,y\in\p$ such that $\hat{\phi}_2(x)=g^{-1}(c)$,   $\hat{\phi}_2(y)=g^{-1}(d)$, and $gf'g^{-1}(x)=y$, this projection is an epimorphism.
Let $\gamma=\bar{\phi}_1$ and $\delta=\bar{\phi}_2\circ g^{-1}$. Then $\bar{\psi}\circ\alpha\circ\gamma=\bar{\psi}\circ\beta\circ\delta$.

\end{proof}

\noindent {\bf{Acknowledgments.}} 
I would like to thank S{\l}awomir Solecki for supplying a proof of Lemma \ref{slawek}.
I would like to thank Todor Tsankov for pointing out a mistake in the published version of the paper and suggesting a correct statement of Theorem \ref{abcd}.

\end{document}